\documentclass[leqno,12pt]{amsart}
\usepackage{latexsym}
\usepackage{mathrsfs}
\usepackage{amsmath}
\usepackage{amssymb}
\usepackage[bookmarks]{hyperref}
\usepackage{pstricks,pst-plot}

\font\tenscr=rsfs10 % scaled \magstep1
\font\sevenscr=rsfs7 % scaled \magstep1
\font\fivescr=rsfs5 % scaled \magstep1
\skewchar\tenscr='177 \skewchar\sevenscr='177 \skewchar\fivescr='177
\newfam\scrfam \textfont\scrfam=\tenscr \scriptfont\scrfam=\sevenscr
\scriptscriptfont\scrfam=\fivescr
\def\scr{\fam\scrfam}

\newtheorem{theorem}{Theorem}[section]
\newtheorem{lemma}[theorem]{Lemma}
\newtheorem{corollary}[theorem]{Corollary}
\newtheorem{proposition}[theorem]{Proposition}

\theoremstyle{definition}

\newtheorem{definition}[theorem]{Definition}

\allowdisplaybreaks[1]

\newcommand{\C}{\mathbb{C}}
\def\C{{\mathbf {C}\/}}

\newcommand{\Z}{\mathbb{Z}}

\newcommand{\R}{\mathbf{R}}
\def\R{{\mathbf {R}\/}}

\def\Z{\mathbb Z}

\def\sF{{\scr F}}

\def \ma{{{\mathfrak{M}}_A}}

\def\tA{\widetilde A}

\def\wX{\widehat X}

\newcommand{\bthm}{\begin{theorem}}
\newcommand{\ethm}{\end{theorem}}
\newcommand{\blem}{\begin{lemma}}
\newcommand{\elem}{\end{lemma}}
\newcommand{\bcor}{\begin{corollary}}
\newcommand{\ecor}{\end{corollary}}
\newcommand{\bprop}{\begin{proposition}}
\newcommand{\eprop}{\end{proposition}}
\newcommand{\bdefn}{\begin{definition}}
\newcommand{\edefn}{\end{definition}}
\newcommand{\bpf}{\begin{proof}}
\newcommand{\epf}{\end{proof}}

\newcommand{\bm}{\bibitem}

\newcommand{\bi}{\begin{itemize}}
\newcommand{\ei}{\end{itemize}}

\newcommand{\bc}{\begin{cases}}
\newcommand{\ec}{\end{cases}}

\newcommand{\ba}{\begin{array}}
\newcommand{\ea}{\end{array}}

\newcommand{\be}{\begin{equation}}
\newcommand{\ee}{\end{equation}}

\newcommand{\bea}{\begin{eqnarray}}
\newcommand{\eea}{\end{eqnarray}}

\newcommand{\beaa}{\begin{eqnarray*}}
\newcommand{\eeaa}{\end{eqnarray*}}

\newcommand{\beastar}{\begin{eqnarray*}}
\newcommand{\eeastar}{\end{eqnarray*}}

\font\tenscr=rsfs10 % scaled \magstep1
\font\sevenscr=rsfs7 % scaled \magstep1
\font\fivescr=rsfs5 % scaled \magstep1
\skewchar\tenscr='177 \skewchar\sevenscr='177 \skewchar\fivescr='177
\newfam\scrfam \textfont\scrfam=\tenscr \scriptfont\scrfam=\sevenscr
\scriptscriptfont\scrfam=\fivescr
\def\scr{\fam\scrfam}

\def\scr{\fam\scrfam}

\def \bs {\setminus}

\def\h#1{\widehat {#1}}
\def\hh#1{\widehat {#1} \bs  {#1}}

\def \C {\mathbb C}

\def\R {\mathbb R}

\def \Z {\mathbb Z}

\def\sF{{\scr F}}

\def \ma {\mathfrak{M}_A}  % or \mathscr *

\begin{document}

\title{A hull that contains no discs}
\author{Alexander J. Izzo}
\address{Department of Mathematics and Statistics, Bowling Green State University, Bowling Green, OH 43403}
\email{aizzo@bgsu.edu}

%\keywords{}


\begin{abstract}
It is shown that if $A$ is a unital commutative Banach algebra with a dense set of invertible elements,
then the maximal ideal space of $A$ contains no compact, locally connected, simply coconnected subspace of topological dimension $\geq 2$.
As a consequence, the existence of a compact set in $\C^2$ with a nontrivial polynomial hull that contains no topological discs is obtained.  This strengthens the celebrated result of Stolzenberg from 1963 that there exists a nontrivial polynomial hull that contains no analytic discs, and it answers a question stated in the \hbox{literature} 10 years ago by Dales and Feinstein but considered much earlier. 
\end{abstract}
\maketitle

%\vskip -2.55 true in
%\centerline{\footnotesize\it Dedicated to Andrew Browder} 
%\vskip 2.55 truein

%========================================================================

\section{Introduction}

It was once conjectured that whenever the polynomial hull $\wX$ of a compact set $X$ in $\C^N$ is strictly larger than $X$, the complementary set 
$\wX\setminus X$ must contain an analytic variety of positive dimension, and hence in particular, an analytic disc.  This conjecture was disproved by
Gabriel~Stolzenberg \cite{Stol}.  
Stolzenberg's result immediately suggests the question of whether the set $\wX\setminus X$, when nonempty, must contain a \emph{topological} disc (i.e., a subspace homeomorphic to a disc in the plane).  We will prove that the answer is \emph{no}.  Henceforth, by a \emph{disc} we shall mean a topological space homeomorphic to a disc in the plane.

\bthm\label{nodisc}
There exists a compact set $X$ in $\C^2$ such that the set $\hh X$ is nonempty but $\h X$ contains no discs.
\ethm

This result sheds a new light on the failure of the existence of analytic discs in hulls since it shows that the failure can occur for purely topological, rather than analytical, reasons.

The question of whether the set $\wX\setminus X$, when nonempty, must contain a disc was stated explicitly 10 years ago in the paper \cite{DalesF} of Garth Dales and Joel Feinstein, but it was likely considered over 50 years ago when Stolzenberg gave his example without \emph{analytic} discs.  The issue was mentioned 20 years ago in the paper \cite{Lev} of Norman Levenberg where it was erroneously asserted that the hull without analytic discs constructed by John Wermer in \cite{Wermerpark} contains no (topological) discs.  This assertion may or may not be true; it has not been proven either way.  The set $X$ constructed by Wermer in \cite{Wermerpark} (and also presented in \cite[pp.~211--218]{AW}) is a subset of $\C^2$.  
Wermer proved that for each constant $c\in \C$, the vertical slice $\{(z_1,z_2)\in X: z_1=c\}$ is totally disconnected; thus $X$ contains no vertical disc.  He also proved that there does not exist an open set $\Omega\subset \C$ and a continuous function$f:\Omega\rightarrow \C$ whose graph $\{(z, f(z)): z\in \Omega\}$ is contained in $X$; thus $X$ contains no disc that is a graph over a set in $\C$.  However, a topological disc, as opposed to a \emph{smooth} disc, need not be locally a graph.  In fact even a \emph{smooth} disc, as opposed to an \emph{analytic} disc, need not be locally a graph \emph{over} $\C$.  For instance, the disc $\Delta=\{(x_1, x_2)\in \R^2\subset \C^2: |x_1|^2 +|x_2|^2 \leq 1/2\}$ is not locally a graph \emph{over} $\C$.  Thus it has not been proven even that Wermer's hull contains no \emph{smooth} discs.  Furthermore, using the construction of hulls without analytic discs given in the paper \cite{DuvalL} of Julien Duval and Levenberg, it is immediate that there is a compact set $X\subset \C^2$ such that $\hh X\supset \Delta$ and $\h X$ contains no analytic discs.

We will obtain Theorem~\ref{nodisc} by proving a general statement about the topology of the maximal ideal space of a unital commutative Banach algebra with dense invertible group.  Following Dales and Feinstein, 
we will say that a Banach algebra $A$ has \emph{dense invertibles} if the invertible elements of $A$ are dense in $A$.  Following Stolzenberg, we will say that a compact space $X$ is {\em simply coconnected\/} if the first \v Cech cohomology group $\check H^1(X,\Z)$ vanishes.  By the \emph{topological dimension} of a space, we shall mean the usual Lebesgue covering dimension.

Our main theorem is as follows.

\bthm\label{general}
Let $A$ be a unital commutative Banach algebra with dense invertibles.
Then the maximal ideal space $\ma$ of $A$ contains no compact, locally connected, simply coconnected subspace of topological dimension $\geq 2$.
%greater than or equal to 2.
In particular, $\ma$ contains no discs.
\ethm

As an immediate corollary we have:

\bcor
Suppose $X$ is a compact set in $\C^N$ such that $P(X)$ has dense invertibles.   Then the polynomial hull $\h X$ contains no compact, locally connected, simply coconnected subspace of topological dimension $\geq 2$.
In particular, $\h X$ contains no discs.
\ecor

Since it was proven by Garth Dales and Joel Feinstein \cite[Theorem~2.1]{DalesF} that there 
exists a compact set $X$ in $\C^2$ such that $\hh X$ is nonempty and $P(X)$ has dense invertible, Theorem~\ref{nodisc} follows as an immediate corollary.
Numerous additional examples of \lq\lq hulls with dense invertibles\rq\rq\ in $\C^N$, having various additional properties, have been given in the papers
\cite{Izzo0} and \cite{Izzo} of the author and in the paper \cite{CGI} of Brian Cole, Swarup Ghosh, and the author.  Thus we immediately get the existence of hulls without discs having these additional properties.  For instance, there exists a Cantor set in $\C^3$ with nontrivial polynomial hull containing no discs; there exists a compact set $X$ in $\C^3$ with nontrivial polynomial hull containing no discs such that $P(X)$ has no nontrivial Gleason parts; there exists a compact set $X$ in $\C^N$ ($N\geq 2$) with nontrivial polynomial hull containing no discs such that $P(X)$ has a \lq\lq large\rq\rq\ Gleason part.

The next result is a variation of Theorem~\ref{general} in which a stronger conclusion is obtained from a stronger hypothesis.  

\bthm\label{special}
Let $A$ be a unital commutative Banach algebra with dense invertibles and such that the set of squares of the elements of $A$ are dense in $A$.  Then the maximal ideal space of $A$ contains no compact, locally connected, subspace of topological dimension $\geq 2$. 
\ethm

Surprisingly, while Theorem~\ref{general} makes no reference to square roots, and examples of nontrivial uniform algebras with dense invertibles can be constructed without the use of square roots, the proof of Theorem~\ref{general} will use Brian Cole's system of square root extensions \cite{Cole} to essentially reduce Theorem~\ref{general} to Theorem~\ref{special}.  While Cole's construction has been used many times to produce examples with particular properties, as far as the author knows, this is the first time the construction has been used to prove a result about algebras already obtained without the use of the construction.

This paper owes a great deal to the work of Dales, Thomas Dawson, and Feinstein.  The author's inspiration to consider Banach algebras with dense invertibles came from the work of Dales and Feinstein in \cite{DalesF}, while the main idea on which the proofs rest was inspired by the proof of Dawson and Feinstein of \cite[Theorem~1.7]{DawsonF}.

\section{The proofs}

\bpf[Proof of Theorem~\ref{special}]
Let $\tA$ denote the uniform closure of the set $\h A$ of Gelfand transforms of the elements of $A$.  Then $\tA$ is a uniform algebra with dense invertibles and such that the set of squares of the elements of $\tA$ are dense in $\tA$.  Also, the maximal ideal space of $\tA$ is homeomorphic to the maximal ideal space of $A$.  Therefore, we may assume without loss of generality that $A$ is a uniform algebra of functions defined on its maximal ideal space $\ma$.

Suppose that $E$ is a compact, locally connected subspace of $\ma$.  Let $B$ be the closure of $A$ restricted to $E$.  Then the set of squares of the elements of $B$ are dense in $B$. Therefore, by a theorem of \v Cirka (see \cite[pp.~131--136]{Stout}) $B=C(E)$.  But $B$, and hence $C(E)$, has dense invertibles.  Consequently, the topological dimension of $E$ is at most~1 by \cite[Proposition 3.3.2]{Pears}.
\epf

\bpf[Proof of Theorem~\ref{general}]
Just as in the proof of Theorem~\ref{special}, we may assume without loss of generality that $A$ is a uniform algebra of functions defined on its maximal ideal space $\ma$.

Set 
$X_0=\ma$ and $A_0=A$.  Following the iterative procedure used by Cole \cite[Theorem~2.5]{Cole} (and also described in \cite[p.~ 201]{Stout}) to obtain uniform algebras with only trivial Gleason parts, we will define a sequence of uniform algebras $\{A_m\}_{m=0}^\infty$.  First let $\sF_0$ be the set  of invertible functions in $A_0$.  Let $p:X_0\times \C^{\sF_0}\rightarrow X_0$ and $p_f:X_0\times \C^{\sF_0} \rightarrow \C$ denote the projections given by $p\bigl(x,(y_g)_{g\in \sF_0}\bigr)=x$ and $p_f\bigl(x,(y_g)_{g\in \sF_0}\bigr)=y_f$.  Define $X_1\subset X_0\times \C^{\sF_0}$ by
$$X_1=\{z\in X_0\times \C^{\sF_0}: p_f^2(z)= f(p(z)) \hbox{\ for all } f\in \sF_0\},$$
and let $A_1$ be the uniform algebra on $X_1$ generated by 
the functions $\{p_f: f\in \sF_0\}$.
Note that since $p_f^2=f\circ p$ on $X_1$, each function of the form $f\circ p$ for $f\in\sF_0$ belong to $A_1$.  Thus by identifying each function $f\in A_0$ with $f\circ p$, we can regard $A_0$ as a subalgebra of $A_1$.  
%Then $p_n^2=f_{0,n}$, so we will denote $p_n$ by $\sqrt{f_{0,n}}$.  
By 
\cite[Theorem~2.1]{DawsonF} the algebra $A_1$ has dense invertibles.  Let $\sF_1$ be the set of invertible functions in $A_1$, and carry out the procedure just preformed again, this time using the uniform algebra $A_1$ and the set of functions $\sF_1$, to get another uniform algebra $A_2$ on a space $X_2$.
Iterating the construction we obtain a sequence of uniform algebras $\{A_m\}_{m=0}^\infty$ on compact spaces $\{X_m\}_{m=0}^\infty$ such that for each $m$ the set $\sF_m$ of invertible elements of $A_m$ is dense in $A_m$.  Each $A_m$ can be regarded as a subalgebra of $A_{m+1}$.  

We now take the direct limit of the system of uniform algebras $\{A_m\}$ to obtain a uniform algebra $A_\omega$ on some compact space $X_\omega$.  If we regard each $A_n$ as a subset of $A_\omega$ in the natural way, and set $\sF=\bigcup \sF_m$, then $\sF$ is a dense set of invertibles in $A_\omega$, and every member of $\sF$ has a square root in $\sF$.  By \cite[Theorem~2.5]{Cole}, the maximal ideal space of $A_\omega$ is $X_\omega$.  Therefore, by Theorem~\ref{special} $X_\omega$ contains no compact, locally connected subspace of topological dimension greater than or equal to 2.  Thus to conclude the proof, it suffices to show that for each simply coconnected subspace $E$ of $\ma=X_0$ there is a subspace of $X_\omega$ homeomorphic to $E$.

Let $E_0$ be a simply coconnected subset of $X_0$.  Then each function $g\in \sF_0$ has a continuous logarithm on $E_0$ and hence, a continuous square root $s_g$ on $E_0$.  Then the map $m:E_0\rightarrow X_0\times \C^{\sF_0}$ given by $m(x)=(x, (s_g(x))_{g\in \sF_0})$ has range contained in $X_1$.  Furthermore, $m$ is obviously continuous and injective, and hence is an embedding since $E_0$ is compact.  Thus $m(E_0)$ is a homeomorphic copy of $E_0$ in $X_1$.  For each $m=1, 2, \ldots$, let $\tau_m:X_m\rightarrow X_{m-1}$ be the obvious project.  Then by induction we have that for each $m$ there is a set $E_m$ in $X_m$ homeomorphic to $E_0$ such that $\tau_m$ takes $E_m$ homeomorphically onto $E_{m-1}$.

We have the inverse system of compact spaces
\[ \cdots \stackrel{\tau_{m+2}}{\longrightarrow} X_{m+1} \stackrel{\tau_{m+1}}{\longrightarrow} X_{m} \stackrel{\tau_{m}}\longrightarrow \cdots \stackrel
{\tau_2}{\longrightarrow} X_1
\stackrel{\tau_1}{\longrightarrow} X_0,
\]
where 
$$X_{m+1}=\{z\in X_m\times \C^{\sF_m}: p_f^2(z)= f(p(z)) \hbox{\ for all } f\in \sF_m\},$$
and $\tau_{m+1}\bigl(x,(y_f)_{f\in \sF_m}\bigr)=x$.  
The space $X_\omega$ is the inverse limit of this inverse system:
$$X_\omega=\{(z_k)_{k=0}^\infty\in \prod_{k=0}^\infty X_k: \tau_{m+1}(z_{m+1})=z_m  \hbox{\ for all } m=0, 1, 2, \ldots\}.$$
Set $E_\omega=X_\omega\cap \prod_{k=0}^\infty E_k$.  The set $E_\omega$ is homeomorphic to $E_0$, for the map $E_\omega\rightarrow E_0$ given by $(z_k)_{k=0}^\infty \mapsto z_0$ is obviously a continuous bijection, and hence a homeomorphism since $E_\omega$ is compact.
\epf

\end{document}